\documentclass{sigcomm-alternate}
\usepackage{amsfonts,latexsym,amsmath,amstext,amssymb,verbatim,epsfig,psfrag}
\usepackage{colordvi,pstcol}
\usepackage{graphicx,psboxit}
\usepackage{psfrag}

\begin{document}
\title{Revisiting the D-iteration method: runtime comparison}

\numberofauthors{3}
\author{
   \alignauthor Dohy Hong\\
   \affaddr{Alcatel-Lucent Bell Labs}\\
   \affaddr{Route de Villejust}\\
   \affaddr{91620 Nozay, France}\\
   \email{\normalsize dohy.hong@alcatel-lucent.com}
   \alignauthor G\'erard Burnside\\
   \affaddr{Alcatel-Lucent Bell Labs}\\
   \affaddr{Route de Villejust}\\
   \affaddr{91620 Nozay, France}\\
   \email{\normalsize gerard.burnside@alcatel-lucent.com}
   \alignauthor Philippe Raoult\\
   \affaddr{4 rue Fabert}\\
   \affaddr{75007 Paris, France}\\
   \email{\normalsize philippe@chandra-conseil.fr}
}

\date{\today}
\maketitle

\begin{abstract}
In this paper, we revisit the D-iteration algorithm in order to better explain different performance results that were observed for the numerical computation of the eigenvector associated to the PageRank score. We revisit here the practical computation cost based on the execution runtime compared to the theoretical number of iterations.
\end{abstract}
\category{G.1.3}{Mathematics of Computing}{Numerical Analysis}[Numerical Linear Algebra]
\category{G.2.2}{Discrete Mathematics}{Graph Theory}[Graph algorithms]
\terms{Algorithms, Performance}
\keywords{Numerical computation; Iteration; Fixed point; Gauss-Seidel; Eigenvector.}
\begin{psfrags}
\section{Introduction}\label{sec:intro}
In this paper, we assume that the readers are already familiar with the idea
of the fluid diffusion associated to the D-iteration \cite{d-algo}
to solve the equation:
$$
X = P.X + B
$$
and its application to PageRank equation \cite{dohy}.

For the general description of alternative or existing iteration methods, one
may refer to \cite{Golub1996, Saad}.

This paper investigates further the analysis done in \cite{revisit}. The main results is that
using the container {\em vector} for the iterators, we obtain better and more predictable performance.

\section{Analysis of the computation cost}\label{sec:cs}
\subsection{C++ Programming environment}
For the evaluation of the computation cost, we used
Linux (Ubuntu) machines:

\verb+Intel(R) Core(TM)2 CPU, U7600, 1.20GHz+, cache size 2048 KB (Linux1, $g++-4.4$)

and

\verb+Intel(R) Core(TM) i5 CPU, M560, 2.67GHz+, cache size 3072 KB (Linux2, $g++-4.6$).

The runtime has been measured based on the library {\em time.h} with the function
$clock()$ (with a precision of 10 ms). The runtime below measures the computation time
from the time we start iterations (time to build the iterators are counted separately
as Initialization time). Note that the initialization time has not been optimized here.
Using binary input format, we observed a gain factor of more than 5 compared to the loading time
that are shown in this paper.

The iterators have been built based on the class {\em vector} from STD library:
\begin{verbatim}
  vector<int> l_out[N];
  vector<int> l_in[N];
\end{verbatim}

Compared to results presented in \cite{revisit} where the class {\em list} was used,
we realized that the use of {\em vector} class was in fact appropriate for the
iteration schemes we use in the context of PageRank equations.
The reason is indeed obvious for programmers: {\em vector} is meant to be used for
variable size vector optimizing the access time to the iterator's value (no pointer
required as for {\em list}).
The results in \cite{revisit} were mainly biased by the fact that the {\em list} has
been built naturally column by column and not per row, because of the input file structure:

\begin{verbatim}
#origin_node  destination_node
	0	993508
	1	999978
	2	999978
	3	999978
	5	4
	6	4
	6	962147
	...
\end{verbatim}

\subsection{Java Programming environment}
For the evaluation of the computation cost in Java, we used only the
Linux2 (Ubuntu) machine and JDK version:
\begin{verbatim}
java version "1.6.0_23"
OpenJDK Runtime Environment (IcedTea6 1.11pre)
                 (6b23~pre11-0ubuntu1.11.10.2)
OpenJDK 64-Bit Server VM (build 20.0-b11, mixed mode)
\end{verbatim}
The runtime has been measured by calling the method:\newline
\verb_getCurrentThreadCpuTime()_ (for longer 
computations, we checked that the value was very close to the one measured by \verb_System.currentTimeMillis()_ which was enough precision compared to the C++ measurements).
We used the same rules to start and stop measuring as was done in the C++ implementation.

For performance reasons we used arrays of primitive \verb_int_ {\em TIntArrayList} from \cite{trove} rather than a classic collection of \verb_Integer_ objects (\verb_ArrayList<ArrayList<Integer>>_):
\begin{verbatim}
  ArrayList<TIntArrayList> l_out;
  ArrayList<TIntArrayList> l_in;
\end{verbatim}

As opposed to the C++ implementation, the graph was read directly from a WebGraph compressed file (see \cite{BoVWFI} and \cite{BRSLLP}), we actually used the Java code to generate the text file parsed by the C++ implementation.

\subsection{Algorithms for evaluation}
The algorithms that we evaluated are:
\begin{itemize}
\item PI: Power iteration (equivalent to Jacobi iteration), using row vectors;
\item PI': Power iteration (equivalent to Jacobi iteration), using column vectors;
\item GS: Gauss-Seidel iteration (cyclic sequence);
\item DI-CYC: D-iteration with cyclic sequence (a node $i$ is selected, if $(F_n)_i>0$);
\item DI-SOP (sub-optimal compromise solution): D-iteration with node selection,
   if $(F_n)_i > r_{n'}\times \#out_i/L$, where $\#out_i$ is the out-degree of $i$
   and $r_{n'}$ is computed per cycle $n'$.
\end{itemize}

\subsection{Notations}
\begin{itemize}
\item L: number of non-null entries (links) of $P$;
\item D: number of dangling nodes (0 out-degree nodes);
\item E: number of 0 in-degree nodes: the 0 in-degree nodes are defined recursively:
  a node $i$, having incoming links from nodes that are all 0 in-degree nodes, is
  also a 0 in-degree node; from the diffusion point of view, those nodes are those
  who converged exactly in finite steps;
\item O: number of loop nodes ($p_{ii} \neq 0$);
\item $\max_{in} = \max_i \#in_i$ (maximum in-degree, the in-degree of $i$ is the number of
  non-null entries of the $i$-th row vector of $P$);
\item $\max_{out} = \max_i \#out_i$ (maximum out-degree, the out-degree of $i$ is the number of
  non-null entries of the $i$-th column vector of $P$).
\end{itemize}

\subsection{Pseudo-codes}
The target error we considered here is $1/N$.
\begin{verbatim}
out[i] := out-degree of node i;
in[i]  := in-degree of node i;
l_out[i] := iterator for column i;
l_in[i]  := iterator for row i.
\end{verbatim}

\subsubsection{Power iteration per row: PI}

\begin{verbatim}
for (int i = 0; i < N; i++){
  x_old[i] = 1.0/N;
}
Loop:
while ( error > target_error ){
  for (int i = 0; i < N; i++){
    x_new[i] = (1-d)/N;
  }
  for (int i = 0; i < N; i++){
    for (vector<int>::iterator j=l_in[i].begin(); 
                          j!=l_in[i].end(); j++){
      x_new[i] += d * x_old[*j]/out[*j];
    }
  }
  error = 0.0;
  for (int i = 0; i < N; i++){
    error += x_new[i];
  }
  for (int i = 0; i < N; i++){
    x_new[i] += (1.0-error)/N;
  }
  error = 0.0;
  for (int i = 0; i < N; i++){
    error += abs(x_new[i]-x_old[i]);
  }
  error *= d/(1-d);
}
\end{verbatim}

\subsubsection{Power iteration per column: PI'}
\begin{verbatim}
for (int i = 0; i < N; i++){
  x_old[i] = 1.0/N;
}
Loop:
while ( error > target_error*(1-d)/d ){
  for (int i = 0; i < N; i++){
    x_new[i] = (1-d)/N;
  }
  for (int i = 0; i < N; i++){
    transit = d * x_old[i]/out[i];
    for (vector<int>::iterator j=l_out[i].begin(); 
                          j!=l_out[i].end(); j++){
      x_new[*j] += transit;
    }
  }
  error = 0.0;
  for (int i = 0; i < N; i++){
    error += x_new[i];
  }
  for (int i = 0; i < N; i++){
    x_new[i] += (1.0-error)/N;
  }
  error = 0.0;
  for (int i = 0; i < N; i++){
    error += abs(x_new[i]-x_old[i]);
  }
  error *= d/(1-d);
}
\end{verbatim}

\subsubsection{Gauss-Seidel: GS}

\begin{verbatim}
for (int i = 0; i < N; i++){
  x[i] = (1-d)/N;
}
while ( error > target_error ){
  error = 0.0;
  for (int i = 0; i < N; i++){
    previous = x[i];
    x[i] = (1-d)/N;
    diag = 1.0;
    for (vector<int>::iterator j=l_in[i].begin(); 
                          j!=l_in[i].end(); j++){
      if ( *j != i ){
        x[i] += d * x[*j]/out[*j];
      } else {
        diag -= d/out[i];
      }
    }
    x[i] /= diag;
    error += x[i] - previous;
  }
  e = 0.0;
  for (int i = 0; i < N; i++){
    if ( out[i] == 0 )
      e += x[i];
    }
  }
  error = error*d/(1 - d - d*e);
}
\end{verbatim}

\subsubsection{D-iteration by cycle: DI-CYC}
\begin{verbatim}
for (int i = 0; i < N; i++){
  hist[i] = 0.0;
  fluid[i] = (1-d)/N;
}
e = 0.0;
while ( error > target_error ){
  for (int i=0; i<N; i++){
    if ( fluid[i] > 0 ){
      if ( loop[i] == 1 ){
        transit = fluid[i]*out[i]/(out[i]-d);
      } else {
        transit = fluid[i];
      }
      hist[i] += transit;
      fluid[i] = 0.0;
      if ( outgoing[i] == 0 )
        e += transit;
        double sent = transit*d/out[i];
        for (vector<int>::iterator j=l_out[i].begin(); 
                             j!=l_out[i].end(); j++){
          if ( *j != i ){
            fluid[*j] += sent;
          }
        }
      }
    }
  }
  error = 0.0;
  for (int i=0; i < N; i++){
    error += fluid[i];
  }
  error /= (1 - d - d*e);
}
\end{verbatim}

\subsubsection{D-iteration based on the average diffusion cost: DI-SOP}

Same as for DI-CYC, replacing the condition:
\begin{verbatim}
  if ( fluid[i] > 0 )
\end{verbatim}
by
\begin{verbatim}
  r = 0.0;
  for (int i=0; i < N; i++){
    r += fluid[i];
  }
  if ( fluid[i] > r/L*out[i] )
\end{verbatim}

\subsection{Dataset 1}
In this section, we use the 
web graph imported from the dataset \verb+uk-2007-05@1000000+
(available on \cite{webgraphit}) which has
41,247,159 links on $10^6$ nodes.

Below we vary $N$ from $10^3$ to $10^6$ extracting from the dataset the
information on the first $N$ nodes.
Few graph properties are summarized in Table \ref{tab:1}.

\begin{table}
\begin{center}
\begin{tabular}{|l|cccccc|}
\hline
N        & L/N  & D/N   & E/N   & O/N & $\max_{in}$ & $\max_{out}$\\
\hline
$10^3$   & 12.9 & 0.041 & 0.032 & 0.236 & 716   & 130\\
$10^4$   & 12.5 & 0.008 & 0.145 & 0.114 & 7982  & 751\\
$10^5$   & 31.4 & 0.027 & 0.016 & 0.175 & 34764 & 3782\\
$10^6$   & 41.2 & 0.046 & 0     & 0.204 & 403441& 4655\\
\hline
\end{tabular}\caption{Extracted graph: $N=10^3$ to $10^6$.}\label{tab:1}
\end{center}
\end{table}

In Table \ref{tab:compaP} and \ref{tab:compaP2} we present the results obtained with Linux1
and Linux2: 
\begin{itemize}
\item the prediction by the number of iterations is quite good for GS;
\item the prediction by the number of iterations is quite good for DI-CYC and DI-SOP
  when the compiler optimization is not used;
\item PI' is much better than GS with compiler optimization;
\item PI' and GS are close without compiler optimization;
\item the compiler optimization can bring a speed-up factor (time2/time1) 4-15;
  the gain factor (9-15 for Linux1, 6-17 for Linux2) for column-vector based methods 
  (PI', DI-CYC, DI-SOP) is more important than
  the gain (4 for Linux1, 5 for Linux2) for row-vector based methods (PI, GS).
\end{itemize}

Using Java, we obtained similar results (Table \ref{tab:compaJavaP2}).

\begin{table}
\begin{center}
\begin{tabular}{|l|ccccc|}
\hline
          & PI  & PI'   & GS      & DI-CYC  & DI-SOP\\
\hline
\hline
 \multicolumn{6}{|l|}{$N=10^3$. Init: 0.05s} \\
\hline
 nb iter  & 28   & 28    & 18.7    & 17.5    & 11.1 \\
 speed-up & 1    & 1.0   & 1.5     & 1.6     & 2.5  \\
\hline
 time1 (s)& 0.02 & 0.01  & 0.02    & 0.01    & 0.00 \\
\hline
 time2 (s)& 0.05 & 0.03  & 0.05    & 0.03    & 0.02 \\
\hline
\hline
 \multicolumn{6}{|l|}{$N=10^4$. Init: 0.2s} \\
\hline
 nb iter & 43    & 43    & 30.7    & 26.4    & 12.0 \\
 speed-up& 1     & 1.0   & 1.4     & 1.6     & 3.6 \\
\hline
 time1 (s)& 0.15 & 0.04  & 0.12    & 0.06    & 0.02 \\
 speed-up & 1    & 3.8   & 1.3     & 2.5     & 7.5 \\
\hline
 time2 (s)& 0.64 & 0.43  & 0.52    & 0.35    & 0.16 \\
 speed-up & 1    & 1.5   & 1.2     & 1.8     & 4.0 \\
\hline
 time2/time1 & 4    & 11   & 4     & 6       & 8 \\
\hline
\hline
 \multicolumn{6}{|l|}{$N=10^5$. Init: 3s} \\
\hline
 nb iter  & 52   & 52    & 36.8    & 34.7    & 14.3 \\
 speed-up & 1    & 1.0   & 1.4     & 1.5     & 3.6 \\
\hline
 time1 (s)& 4.5  & 0.83  & 3.5     & 1.1     & 0.52 \\
 speed-up & 1    & 5.4   & 1.3     & 4.1     & 8.7 \\
\hline
 time2 (s)& 18.0 & 12.1  & 14.1    & 9.8     & 4.2 \\
 speed-up & 1    & 1.5   & 1.3     & 1.8     & 4.3 \\
\hline
 time2/time1 & 4    & 15   & 4     & 9       & 8 \\
\hline
\hline
 \multicolumn{6}{|l|}{$N=10^6$. Init: 31s} \\
\hline
 nb iter  & 66   & 66    & 41.8    & 39.8    & 14.6 \\
 speed-up & 1    & 1.0   & 1.6     & 1.7     & 4.5 \\
\hline
 time1 (s)& 75   & 13    & 51      & 16      & 6.3 \\
 speed-up & 1    & 5.8   & 1.5     & 4.7     & 11.9 \\
\hline
 time2 (s)& 296  & 199   & 207     & 144     & 54 \\
 speed-up & 1    & 1.5   & 1.4     & 2.1     & 5.5 \\
\hline
 time2/time1 & 4    & 15   & 4     & 9       & 9 \\
\hline
\end{tabular}\caption{Linux1: Comparison of the runtime for a target error of $1/N$. Speed-up: gain factor w.r.t. PI. time1: with compiler optimization. time2: no compiler optimization.}\label{tab:compaP}
\end{center}
\end{table}

\begin{table}
\begin{center}
\begin{tabular}{|l|ccccc|}
\hline
          & PI  & PI'   & GS      & DI-CYC  & DI-SOP\\
\hline
\hline
 \multicolumn{6}{|l|}{$N=10^3$. Init: 0.01s} \\
\hline
 nb iter  & 28   & 28    & 18.7    & 17.5    & 11.1 \\
 speed-up & 1    & 1.0   & 1.5     & 1.6     & 2.5  \\
\hline
 time1 (s)& 0.01 & 0.00  & 0.00    & 0.00    & 0.00 \\
\hline
 time2 (s)& 0.01 & 0.01  & 0.01    & 0.01    & 0.01 \\
\hline
\hline
 \multicolumn{6}{|l|}{$N=10^4$. Init: 0.05s} \\
\hline
 nb iter & 43    & 43    & 30.7    & 26.4    & 12.0 \\
 speed-up& 1     & 1.0   & 1.4     & 1.6     & 3.6 \\
\hline
 time1 (s)& 0.04 & 0.01  & 0.03    & 0.02    & 0.01 \\
 speed-up & 1    & 4.0   & 1.3     & 2.0     & 4.0 \\
\hline
 time2 (s)& 0.20 & 0.17  & 0.17    & 0.12    & 0.06 \\
 speed-up & 1    & 1.2   & 1.2     & 1.7     & 3.3 \\
\hline
 time2/time1 & 5    & 17   & 6     & 6       & 6 \\
\hline
\hline
 \multicolumn{6}{|l|}{$N=10^5$. Init: 1s} \\
\hline
 nb iter  & 52   & 52    & 36.8    & 34.7    & 14.3 \\
 speed-up & 1    & 1.0   & 1.4     & 1.5     & 3.6 \\
\hline
 time1 (s)& 1.3  & 0.32  & 0.98    & 0.46    & 0.23 \\
 speed-up & 1    & 4.1   & 1.3     & 2.8     & 5.7 \\
\hline
 time2 (s)& 6.1  & 4.7   & 4.5     & 3.5     & 1.6 \\
 speed-up & 1    & 1.3   & 1.4     & 1.7     & 3.8 \\
\hline
 time2/time1 & 5    & 15   & 5     & 8       & 7 \\
\hline
\hline
 \multicolumn{6}{|l|}{$N=10^6$. Init: 11s} \\
\hline
 nb iter  & 66   & 66    & 41.8    & 39.8    & 14.6 \\
 speed-up & 1    & 1.0   & 1.6     & 1.7     & 4.5 \\
\hline
 time1 (s)& 21   & 5.6   & 14      & 6.4     & 3.1 \\
 speed-up & 1    & 3.8   & 1.5     & 3.3     & 6.8 \\
\hline
 time2 (s)& 97   & 78    & 67      & 53      & 21 \\
 speed-up & 1    & 1.2   & 1.4     & 1.8     & 4.6 \\
\hline
 time2/time1 & 5    & 14   & 5     & 8       & 7 \\
\hline
\end{tabular}\caption{Linux2: Comparison of the runtime for a target error of $1/N$. Speed-up: gain factor w.r.t. PI. time1: with compiler optimization. time2: no compiler optimization.}\label{tab:compaP2}
\end{center}
\end{table}

\begin{table}
\begin{center}
\begin{tabular}{|l|ccccc|}
\hline
          & PI  & PI'   & GS      & DI-CYC  & DI-SOP\\
\hline
\hline
 \multicolumn{6}{|l|}{$N=10^3$. Init: 1.9s} \\
\hline
 nb iter  & 28   & 28    & 18.7    & 17.5    & 11.1 \\
 speed-up & 1    & 1.0   & 1.5     & 1.6     & 2.5  \\
\hline
 time (s)& 0.00 & 0.00  & 0.00    & 0.00    & 0.00 \\
\hline
\hline
 \multicolumn{6}{|l|}{$N=10^4$. Init: 2.0s} \\
\hline
 nb iter & 43    & 43    & 30.7    & 26.4    & 11.8 \\
 speed-up& 1     & 1.0   & 1.4     & 1.6     & 3.6 \\
\hline
 time (s)& 0.04 & 0.03  & 0.04    & 0.03    & 0.03 \\
 speed-up & 1    & 1.3   & 1       & 1.3     & 1.3 \\
\hline
\hline
 \multicolumn{6}{|l|}{$N=10^5$. Init: 2.1s} \\
\hline
 nb iter  & 52   & 52    & 36.8    & 34.7    & 14.4 \\
 speed-up & 1    & 1.0   & 1.4     & 1.5     & 3.6 \\
\hline
 time (s)& 1.26  & 0.45  & 0.91    & 0.43    & 0.33 \\
 speed-up & 1    & 2.8   & 1.4     & 2.9     & 3.8 \\
\hline
\hline
 \multicolumn{6}{|l|}{$N=10^6$. Init: 2.8s} \\
\hline
 nb iter  & 66   & 66    & 41.8    & 39.8    & 14.6 \\
 speed-up & 1    & 1.0   & 1.6     & 1.7     & 4.5 \\
\hline
 time (s)& 21.3  & 7.85  & 13.33 & 6.2     & 4.0 \\
 speed-up & 1    & 2.7   & 1.6     & 3.4     & 5.3 \\
\hline
\end{tabular}\caption{Linux2 in Java: Comparison of the runtime for a target error of $1/N$. Speed-up: gain factor w.r.t. PI. Similar to table \ref{tab:compaP2} but in Java!}\label{tab:compaJavaP2}
\end{center}
\end{table}

\subsection{Dataset 1bis}
We considered here the same dataset than dataset 1 but for $P^t$ (transposed matrix, which
means we inverse incoming and outgoing links).
In Table \ref{tab:compaPt} and \ref{tab:compaPt2} we present the results obtained for $P$ with Linux1
and Linux2: 
\begin{itemize}
\item the prediction by the number of iterations is not bad for GS;
\item the prediction by the number of iterations is still good for DI-CYC and DI-SOP
  when the compiler optimization is not used;
\item PI' is still much better than GS with compiler optimization;
\item PI' and GS are still close without compiler optimization;
\item the compiler optimization can bring a speed-up factor (time2/time1) 4-16;
  the gain factor for column-vector based methods 
  (PI': 11-16, DI-CYC and DI-SOP: 5-9) is more important than
  the gain for row-vector based methods (PI and GS: 4-5).
\item the gain factors are globally more stable than we expected compared to results for $P$
  (we expected worse results): this suggests that the performance of the D-iteration approaches
  are quite stable w.r.t. the variance of in-degree/out-degree.
\end{itemize}

The results with Java are shown in Table \ref{tab:compaJavaPt2}.
In the Java implementation the init phase was almost constant (from 1.9s with $N=10^3$ to 2.9s with $N=10^6$) 
when varying $N$, because we always read the full graph file ($10^6$).

\begin{table}
\begin{center}
\begin{tabular}{|l|ccccc|}
\hline
         & PI  & PI'   & GS      & DI-CYC  & DI-SOP\\
\hline
\hline
 \multicolumn{6}{|l|}{$N=10^3$. Init: 0.05s} \\
\hline
 nb iter  & 30   & 30    & 22.6    & 14.4    & 12.2 \\
 speed-up & 1    & 1.0   & 1.3     & 2.1     & 2.5 \\
\hline
 time1 (s)& 0.02 & 0.01  & 0.01     & 0.01     & 0.00 \\
\hline
 time2 (s)& 0.07 & 0.03  & 0.05     & 0.04     & 0.02 \\
\hline
\hline
 \multicolumn{6}{|l|}{$N=10^4$. Init: 0.2s} \\
\hline
 nb iter  & 40   & 40    & 30.7    & 28.3    & 11.5 \\
 speed-up & 1    & 1.0   & 1.3     & 1.4     & 3.5 \\
\hline
 time1 (s)& 0.15 & 0.03  & 0.13    & 0.06    & 0.02 \\
 speed-up & 1    & 5.0   & 1.2     & 2.5     & 7.5 \\
\hline
 time2 (s)& 0.60 & 0.41  & 0.52    & 0.38    & 0.17 \\
 speed-up & 1    & 1.5   & 1.2     & 1.6     & 3.5 \\
\hline
 time2/time1 & 4 & 14    & 4       & 6       & 9 \\
\hline
\hline
 \multicolumn{6}{|l|}{$N=10^5$. Init: 3s} \\
 nb iter  & 51   & 51    & 37.8    & 35.2    & 17.1 \\
 speed-up & 1    & 1     & 1.3     & 1.4     & 3.0 \\
\hline
 time1 (s)& 4.4  & 0.82  & 3.5     & 1.2     & 0.6 \\
 speed-up & 1    & 5.4   & 1.3     & 3.7     & 7.3 \\
\hline
 time2 (s)& 17.3 & 11.9  & 14.2    & 10.4    & 5.0 \\
 speed-up & 1    & 1.5   & 1.2     & 1.7     & 3.5 \\
\hline
 time2/time1 & 4 & 15    & 4       & 9       & 8 \\
\hline
\hline
 \multicolumn{6}{|l|}{$N=10^6$. Init: 31s} \\
\hline
 nb iter  & 78   & 78    & 45.8    & 43.1    & 17.1 \\
 speed-up & 1    & 1.0   & 1.7     & 1.8     & 4.6 \\
\hline
 time1 (s)& 88   & 17.8  & 54.7    & 18.7    & 8.1 \\
 speed-up & 1    & 4.9   & 1.6     & 4.7     & 11 \\
\hline
 time2 (s)& 346  & 238   & 221     & 156     & 65 \\
 speed-up & 1    & 1.5   & 1.6     & 2.2     & 5.3 \\
\hline
 time2/time1 & 4 & 13    & 4       & 8       & 8 \\
\hline
\end{tabular}\caption{$P^t$ on Linux1: Comparison of the runtime for a target error of $1/N$. Speed-up: gain factor w.r.t. PI. time1: with compiler optimization. time2: no compiler optimization.}\label{tab:compaPt}
\end{center}
\end{table}

\begin{table}
\begin{center}
\begin{tabular}{|l|ccccc|}
\hline
         & PI  & PI'   & GS      & DI-CYC  & DI-SOP\\
\hline
\hline
 \multicolumn{6}{|l|}{$N=10^3$. Init: 0.01s} \\
\hline
 nb iter  & 30   & 30    & 22.6    & 21.1    & 12.2 \\
 speed-up & 1    & 1.0   & 1.3     & 1.4     & 2.5 \\
\hline
 time1 (s)& 0.00 & 0.00  & 0.00     & 0.00     & 0.00 \\
\hline
 time2 (s)& 0.02 & 0.01  & 0.02     & 0.01     & 0.01 \\
\hline
\hline
 \multicolumn{6}{|l|}{$N=10^4$. Init: 0.05s} \\
\hline
 nb iter  & 40   & 40    & 30.7    & 28.3    & 11.5 \\
 speed-up & 1    & 1.0   & 1.3     & 1.4     & 3.5 \\
\hline
 time1 (s)& 0.04 & 0.01  & 0.04    & 0.03    & 0.01 \\
 speed-up & 1    & 4.0   & 1.0     & 1.3     & 4.0 \\
\hline
 time2 (s)& 0.19 & 0.16  & 0.17    & 0.14    & 0.06 \\
 speed-up & 1    & 1.2   & 1.1     & 1.4     & 3.2 \\
\hline
 time2/time1 & 5 & 16    & 4       & 5       & 6 \\
\hline
\hline
 \multicolumn{6}{|l|}{$N=10^5$. Init: 1s} \\
 nb iter  & 51   & 51    & 37.8    & 35.2    & 17.1 \\
 speed-up & 1    & 1     & 1.3     & 1.4     & 3.0 \\
\hline
 time1 (s)& 1.2  & 0.41  & 0.93    & 0.54    & 0.31 \\
 speed-up & 1    & 2.9   & 1.3     & 2.2     & 3.9 \\
\hline
 time2 (s)& 5.5  & 4.8   & 4.6     & 3.7     & 2.1 \\
 speed-up & 1    & 1.1   & 1.2     & 1.5     & 2.6 \\
\hline
 time2/time1 & 5 & 12    & 5       & 7       & 7 \\
\hline
\hline
 \multicolumn{6}{|l|}{$N=10^6$. Init: 11s} \\
\hline
 nb iter  & 78   & 78    & 45.8    & 43.1    & 17.1 \\
 speed-up & 1    & 1.0   & 1.7     & 1.8     & 4.6 \\
\hline
 time1 (s)& 24   & 8.8   & 15      & 8.4     & 4.3 \\
 speed-up & 1    & 2.7   & 1.6     & 2.9     & 5.6 \\
\hline
 time2 (s)& 111  & 95    & 72      & 58      & 27 \\
 speed-up & 1    & 1.2   & 1.5     & 1.9     & 4.1 \\
\hline
 time2/time1 & 5 & 11    & 5       & 7       & 6 \\
\hline
\end{tabular}\caption{$P^t$ on Linux2: Comparison of the runtime for a target error of $1/N$. Speed-up: gain factor w.r.t. PI. time1: with compiler optimization. time2: no compiler optimization.}\label{tab:compaPt2}
\end{center}
\end{table}

\begin{table}
\begin{center}
\begin{tabular}{|l|ccccc|}
\hline
         & PI  & PI'   & GS      & DI-CYC  & DI-SOP\\
\hline
\hline
 \multicolumn{6}{|l|}{$N=10^3$. Init: 1.9s} \\
\hline
 nb iter  & 30   & 30    & 22.6    & 21.1    & 12.2 \\
 speed-up & 1    & 1.0   & 1.3     & 1.4     & 2.5 \\
\hline
 time (s)& 0.01 & 0.00  & 0.01     & 0.00     & 0.00 \\
\hline
\hline
 \multicolumn{6}{|l|}{$N=10^4$. Init: 2.0s} \\
\hline
 nb iter  & 40   & 40    & 30.7    & 28.3    & 11.4 \\
 speed-up & 1    & 1.0   & 1.3     & 1.4     & 3.5 \\
\hline
 time (s)& 0.05 & 0.02  & 0.04    & 0.03    & 0.02 \\
 speed-up & 1    & 2.5   & 1.25    & 1.7     & 2.5 \\
\hline
\hline
 \multicolumn{6}{|l|}{$N=10^5$. Init: 2.1s} \\
\hline
 nb iter  & 51   & 51    & 37.8    & 35.2    & 17.4 \\
 speed-up & 1    & 1     & 1.3     & 1.4     & 2.9 \\
\hline
 time (s)& 1.23  & 0.50  & 0.92    & 0.46    & 0.38 \\
 speed-up & 1    & 2.5   & 1.3     & 2.7     & 3.2 \\
\hline
\hline
 \multicolumn{6}{|l|}{$N=10^6$. Init: 2.9s} \\
\hline
 nb iter  & 78   & 78    & 45.8    & 43.1    & 17.2 \\
 speed-up & 1    & 1.0   & 1.7     & 1.8     & 4.5 \\
\hline
 time (s)& 25.6   & 9.7   & 15.3   & 7.4     & 5.1 \\
 speed-up & 1    & 2.6   & 1.7     & 3.5     & 5.0 \\
\hline
\end{tabular}\caption{$P^t$ on Linux2 in Java: Comparison of the runtime for a target error of $1/N$. Speed-up: gain factor w.r.t. PI. Similar to table \ref{tab:compaPt2} but in Java!}\label{tab:compaJavaPt2}
\end{center}
\end{table}

\subsection{Dataset 2}\label{sec:calif}

Below, we used the web graph \verb+gr0.California+
(available on \verb+http://www.cs.cornell.edu/Courses/cs685/+ \verb+2002fa/+).
The main motivation was here to try to understand the unexpected (too much) gain
observed in \cite{dohy} for this graph.

\begin{table}
\begin{center}
\begin{tabular}{|cccccc|}
\hline
L/N & D/N & E/N & O/N & $\max_{in}$ & $\max_{out}$\\
\hline
1.67& 0.48& 0.91& 0   & 199         & 164\\
\hline
\end{tabular}\caption{$N=9664$.}\label{tab:cal}
\end{center}
\end{table}

As it has been pointed out in \cite{revisit} (Table \ref{tab:cal}), this graph is very specific in that
more than 90\% of nodes are 0 in-degree nodes.
The runtime is here too short to make a comparison.

\begin{table}
\begin{center}
\begin{tabular}{|r|ccccc|}
\hline
         & PI  & PI'   & GS      & DI-CYC  & DI-SOP\\
\hline
 \multicolumn{6}{|l|}{$N=9664$. Init: 0.1s} \\
\hline
 nb iter & 43  & 43    & 22      & 3.1   & 1.8\\
 speed-up& 1   & 1.0   & 2.0      & 14    & 24\\
\hline
 time1   & 0.03 & 0.02 & 0.04    & 0.01   & 0.01\\
 speed-up& 1   & 1.5   & 0.8     & 3.0    & 3.0\\
\hline
 time2   & 0.16& 0.12  & 0.09    & 0.04  & 0.03\\
 speed-up& 1   & 1.3   & 1.8     & 4.0   & 5.3\\
\hline
 time2/time1& 5 & 6    & 2       & 4     & 3\\
\hline
\end{tabular}\caption{$P$. times1: with compiler optimization, time2: without compiler optimization.}\label{tab:compa-calif}
\end{center}
\end{table}

Table \ref{tab:compa-calif-inv} presents the results of the computation cost
associated to the matrix $P^t$ for comparison.

\begin{table}
\begin{center}
\begin{tabular}{|r|ccccc|}
\hline
         & PI  & PI'   & GS      & DI-CYC  & DI-SOP\\
\hline
 \multicolumn{6}{|l|}{$N=9664$. Init: 0.1s} \\
\hline
 nb iter & 28  & 28    & 16      & 5.6   & 2.0\\
 speed-up& 1   & 1.0   & 1.8     & 5.0   & 14\\
\hline
 time1   & 0.04& 0.01  & 0.03    & 0.01  & 0.01\\
 speed-up& 1   & 4     & 1.3     & 4     & 4\\
\hline
 time2   & 0.11& 0.07  & 0.08    & 0.04  & 0.03\\
 speed-up& 1   & 1.6   & 1.4     & 2.8   & 3.7\\
\hline
 time2/time2& 3& 7     & 3       & 4     & 3\\
\hline
\end{tabular}\caption{$P^t$. times1: with compiler optimization, time2: without compiler optimization.}\label{tab:compa-calif-inv}
\end{center}
\end{table}

\section{Conclusion}\label{sec:conclusion}
In this paper we revisited the D-iteration method with a practical consideration
of the computation cost to solve the PageRank equation for web graphs: 
the use of the class {\em vector} for iterators produced much faster results
with performance that are closer to expectations.

\end{psfrags}
\bibliographystyle{abbrv}
\bibliography{sigproc}

\end{document}